\documentclass{article}
\usepackage[a4paper]{geometry}
\usepackage{setspace}
\usepackage{authblk}
\usepackage[utf8]{inputenc}
\usepackage[T1]{fontenc}
\usepackage{amsmath,amsfonts,amsthm,amssymb}
\usepackage{natbib}
\newenvironment{acks}[1][Acknowledgments]{\section*{#1}}{}
\newenvironment{funding}{\section*{Funding}}{}
\newenvironment{supplement}{\section*{Supplementary Material}}{}
\newcommand{\stitle}[1]{\noindent\textbf{#1}\par}
\newcommand{\sdescription}[1]{\noindent #1\par}

\usepackage{amsmath,amssymb}
\usepackage{microtype}
\usepackage{booktabs,array}
\usepackage{etoolbox}
\AtBeginEnvironment{tabular}{\setlength{\tabcolsep}{3.4pt}}
\usepackage{graphicx}
\usepackage{algorithm}
\usepackage{algpseudocode}
\usepackage{xurl}
\usepackage{hyperref}

\newcommand{\bdoi}[1]{\url{https://doi.org/#1}}

\numberwithin{equation}{section}

\theoremstyle{plain}
\newtheorem{theorem}{Theorem}[section]
\newtheorem{lemma}[theorem]{Lemma}
\newtheorem{proposition}[theorem]{Proposition}
\newtheorem{corollary}[theorem]{Corollary}

\theoremstyle{definition}
\newtheorem{assumption}{Assumption}

\newtheorem{remark}[theorem]{Remark}

\newcommand{\E}{\mathbb E}
\newcommand{\Prob}{\mathbb P}
\newcommand{\one}{\mathbf 1}
\newcommand{\Simp}{\mathcal S_{n,K}}
\newcommand{\Perm}{\mathfrak S_K}
\newcommand{\cB}{\mathcal B}
\newcommand{\cT}{\mathcal T}
\newcommand{\cC}{\mathcal C}

\newcommand{\pistar}{\pi_n^*}
\newcommand{\clip}{\operatorname{clip}}
\newcommand{\diag}{\operatorname{diag}}

\newcommand{\Bern}{\operatorname{Bern}}
\newcommand{\Cat}{\operatorname{Categorical}}

\newcommand{\normop}[1]{\lVert #1\rVert_{\mathrm{op}}}
\newcommand{\normF}[1]{\lVert #1\rVert_{\mathrm F}}
\newcommand{\normmax}[1]{\lVert #1\rVert_{\max}}
\newcommand{\normone}[1]{\lVert #1\rVert_1}
\newcommand{\argmax}{\operatorname*{arg\,max}}

\newcommand{\econst}{\mathrm e}

\begin{document}
\title{Optimal Community Recovery by Spectrally Initialized Variational EM in General Stochastic Block Models}

\author[A]{Binghui~Liu}
\author[B]{Jiangzhou~Wang\thanks{Corresponding author. Email: \texttt{wangjz695@szu.edu.cn}}}
\author[C]{Man-Chung~Yue}
\affil[A]{School of Mathematics and Statistics \& KLAS,
Northeast Normal University, Changchun, China, \texttt{liubh100@nenu.edu.cn}}
\affil[B]{School of Mathematical Sciences, Institute of Statistical Sciences,
Shenzhen University, Shenzhen 518060, China, \texttt{wangjz695@szu.edu.cn}}
\affil[C]{Department of Applied Mathematics,
The Hong Kong Polytechnic University, Hung Hom, Hong Kong, \texttt{manchung.yue@polyu.edu.hk}}
\date{}
\maketitle

\begin{abstract}
We prove a Chernoff-exponent guarantee for the output of spectrally initialized batch variational EM after a prescribed iteration budget, without assuming global maximization of its variational objective. The iteration repeatedly estimates the entire block probability matrix and community proportions from the same sparse graph. We consider a fixed number of communities with proportions bounded away from zero and fixed, positive, distinct connectivity profiles; neither assortativity nor full rank is required. The algorithm uses simultaneous softmax updates without sample splitting or posterior thresholding. Its analysis must control the feedback between estimated parameters, soft labels and reused edges at an exponentially small risk scale. We establish a uniform one-step bound over data-dependent soft assignments whose random remainder has exponentially small expectation. Combined with an exponentially reliable regularized spectral initializer, this yields an unconditional expected misclassification rate bounded by $\exp\{-(1-o(1))J_n\}$ throughout the sparse, diverging-degree regime, where $J_n$ is the minimum nodewise Chernoff information. Matching lower bounds establish first-order logarithmic minimax optimality on local parameter spaces allowing unknown connectivity and varying community counts. The algorithm also attains the sharp first-order exact-recovery threshold. Numerical experiments illustrate refinement gains and sensitivity to initialization and imbalance.
\end{abstract}


\section{Introduction}\label{sec:introduction}
Statistical theory for community-detection models includes guarantees for estimators defined by global optima of likelihood or variational criteria \citep{CelisseEtAl2012,BickelEtAl2013}. Such results concern an idealized optimizer. Finding it can involve combinatorial or nonconvex optimization, and an iterative method need not attain it \citep{BleiEtAl2017}. The optimizer's consistency or statistical accuracy therefore does not automatically describe the estimate returned by that method. This distinction does not imply that the algorithmic estimate is inconsistent; it means that its risk needs a separate analysis. Existing algorithmic theory has already supplied such analyses for iterative variational inference and optimal refinement \citep{ZhangZhou2020,ZhouLi2018}. Our question is whether the prescribed batch variational EM output can attain optimal recovery accuracy without first establishing global optimization.

We give an affirmative answer for $B^*=\rho_nQ$ with fixed $K$, fixed positive symmetric $Q$ having distinct rows, unknown proportions bounded away from zero, and $d_n=n\rho_n\to\infty$, $\rho_n\to0$. Neither assortativity nor full rank is required. Starting from the prescribed regularized spectral initializer, the algorithm repeatedly estimates the full block matrix and proportions and makes simultaneous soft updates on the same graph, without sample splitting or posterior thresholding. After the prescribed iteration budget, its returned estimate satisfies $\E R(\widehat z,z^*)\le\exp\{-(1-o(1))J_n\}$, where $R$ is the permutation-invariant misclassification proportion and $J_n$ is the minimum Chernoff information for classifying one vertex when the other labels and parameters are known.

\subsection{Contributions and proof strategy}
\emph{Uniform control of the full plug-in iteration.} Unlike a planted-partition model, where all community pairs share a single within-/between-block contrast, distinguishing communities \(k\) and \(\ell\) here depends on the full weighted difference between their connectivity profiles.Estimating the full matrix \(B\) perturbs the pair-specific log-odds used in the E-step. At the same time, the current soft assignments used to estimate \(B\) are themselves functions of the same observed edges, creating a feedback between edge noise, parameter estimation, and label updates. To obtain a guarantee for the computed iterate, we control this feedback uniformly over every soft input in a fixed average-loss neighborhood (Theorem~\ref{thm:contraction}). Consequently, the bound applies to all data-dependent iterates computed from the same original graph, without splitting the data or replacing soft memberships by hard labels between updates.

\emph{An unconditional expected-risk exponent.} Uniform contraction on a high-probability event is not enough to obtain the target expected-risk exponent: the event's failure probability may vanish and yet still be much larger than \(e^{-J_n}\). The main technical estimate, Lemma~\ref{lem:onestep}, retains a common random remainder whose expectation is bounded by $Ce^{-J_n+h}$. It combines uniform soft-count control, regularized operator-norm control of the noise--membership interaction, and an oracle likelihood-ratio bound used only in the proof. Together with an exponentially reliable initializer, this remainder allows us to average over the full graph distribution and obtain the unconditional expected-risk guarantee in Theorem~\ref{thm:refinement} throughout the diverging-degree sparse regime.

\emph{A matching benchmark for the specified algorithm.}
On local parameter spaces allowing both $Q$ and community counts to vary,
Theorem~\ref{thm:minimax} identifies the same first-order logarithmic minimax risk.
Theorem~\ref{thm:exact} shows that this iteration also attains the known
first-order exact-recovery threshold
\citep{Abbe2018,AbbeSandon2015limits}.
These comparisons establish first-order sharpness for the specified procedure.
Specifically, this sharpness refers to matching the leading logarithmic error rate
on the stated local parameter spaces and attaining the first-order exact-recovery
threshold; it does not mean that the iteration finds a global maximizer of the
variational objective.
The spectral and variational components are established methods; the contribution
is their sharp statistical analysis at a prescribed iteration budget with full
parameter learning and repeated graph reuse.

\subsection{Comparison with existing theory}
The closest iterative predecessor is \citet{ZhangZhou2020}. Their main planted-partition iteration bounds establish same-graph batch contraction and high-probability minimax-order accuracy while learning the within- and between-community probabilities; it does not require exact balance. Relative to those guarantees, we treat the complete unknown block matrix and estimated proportions under general connectivity, and obtain an unconditional expected Chernoff-exponent bound. This requires controlling pair-specific plug-in log-odds rather than one common block contrast.

The comparison with \citet{ZhouLi2018} is different. They already derive general-SBM Chernoff risk bounds and analyze likelihood-ratio refinement, with different adjacency subblocks used at key stages of their Algorithm~1 to control dependence. We instead analyze repeated simultaneous soft updates on the full graph, with the entire $B$ and the proportions re-estimated at each step. Their treatment obtains matching Chernoff-risk bounds up to multiplicative constants under their assumptions, whereas our sharpness concerns the first-order logarithm. \citet{LiLe2023} give the general soft-count and softmax formulas and study posterior-thresholded batch updates. Our guarantee concerns the unthresholded iteration. 

Consistency and asymptotic parameter theory for likelihood and variational estimators are developed by \citet{CelisseEtAl2012,BickelEtAl2013}, and \citet{GaucherKlopp2021} study variational optimality with missing links. Initialization and parameter learning can alter batch-update behavior even in two-class models \citep{MukherjeeEtAl2018,SarkarEtAl2021}. Exponential minimax theory \citep{ZhangZhou2016}, optimal refinement \citep{GaoEtAl2017}, discrete-structure iteration \citep{GaoZhang2022}, and general-SBM recovery algorithms \citep{AbbeSandon2015limits,AbbeSandon2015agnostic} supply further statistical and computational benchmarks. Our lower bounds show that the specified VEM iterates attain the first-order minimax benchmark on an explicit parameter space that allows local variation in both the connectivity matrix and community counts.

The initializer uses existing reduced-rank geometry: $\widehat U_K\widehat\Lambda_K$ is the symmetric counterpart of SC-RRE in \citet[Algorithm~4 and Theorem~6]{ZhouAmini2019}. Eigenvalue and singular-value weighting have identical rowwise distances, since they differ only by an orthogonal sign matrix. We combine this geometry with degree regularization \citep{LeLevinaVershynin2017} and matrix tools \citep{LeiRinaldo2015,BandeiraVanHandel2016,Meckes2004}, calibrating initialization failure to the expected-risk scale. The coarse spectral bound needed here is not a limit on spectral methods; exponential-error results are available \citep{ZhangSpectral2023}. The connection between likelihood and optimal transport studied by \citet{QuericEtAl2026} addresses complementary questions in optimization and model selection.

The rest of the article is organized as follows. Sections~\ref{sec:model}--\ref{sec:optimality} introduce the model and algorithm and present the main recovery guarantees. Section~\ref{sec:mechanism} develops the key estimates underlying the expected-risk theorem, while Section~\ref{sec:computation} discusses computational aspects and the scope of the theory. Section~\ref{sec:numerics} reports the numerical results, including both refinement gains and settings that reveal practical limitations. All proofs are collected in Sections~A--J of the supplementary material \citep{SVEMSupplement2026}, with additional numerical results and implementation details in Section~K.

\section{Model and Algorithm}\label{sec:model}

This section introduces three ingredients in the order needed below: the network model, the loss and information used to evaluate recovery, and the updates that the algorithm actually computes. Quantities defined from the truth are analytical benchmarks; they are not additional inputs to the algorithm.
\subsection{A general block model with unknown counts}
Let $K\ge2$ be fixed. Conditional on a deterministic label vector $z^*\in[K]^n$, the variables $\{A_{ij}:i<j\}$ are independent and
\begin{equation}\label{eq:model}
 A_{ij}\sim\Bern(B^*_{z_i^*,z_j^*}),\qquad A_{ji}=A_{ij},\qquad A_{ii}=0.
\end{equation}
Write $n_k=\sum_i\one\{z_i^*=k\}$ and $\pi_{k,n}^*=n_k/n$.

\begin{assumption}[Positive community proportions]\label{ass:proportions}
For a fixed $c_0\in(0,1/K)$, $\min_k\pi_{k,n}^*\ge c_0$. The counts $n_k$ are not given to the algorithm and are not imposed as constraints on its output.
\end{assumption}
\begin{assumption}[Positive, separated connectivity profiles]\label{ass:Q}
$B^*=\rho_n Q$, where $Q=Q^\top$ is fixed and all its entries are positive. Its rows are pairwise distinct. Define the actual extrema
\[
 q_-=\min_{k,\ell}q_{k\ell}>0,\qquad q_+=\max_{k,\ell}q_{k\ell}<\infty.
\]
No assortativity or full-rank condition is imposed.
\end{assumption}
\begin{assumption}[Sparse, diverging-degree regime]\label{ass:sparsity}
$\rho_n\to0$ and $d_n=n\rho_n\to\infty$.
\end{assumption}

The scale $d_n$ is proportional to the average expected degree, not necessarily equal to it. The factorization $B^*=\rho_nQ$ is used for asymptotic analysis; the algorithm estimates $B^*$ rather than attempting to identify its scale and shape separately. We write
\begin{equation}\label{eq:P}
\begin{split}
 D_n&=q_+d_n,\\
 P=\E A&=P_0-\diag\{\rho_nq_{z_i^*,z_i^*}\}_{i=1}^n,\\
 P_0&=\rho_n Z^*QZ^{*\top},\qquad W=A-P,
\end{split}
\end{equation}
where $Z^*_{ik}=\one\{z_i^*=k\}$. Probabilities and expectations are conditional on $z^*$ unless otherwise stated. All asymptotic bounds are understood to hold for sufficiently large $n$.

The distinction between $P$ and $P_0$ will matter in the spectral proof. The matrix $P_0$ has rank at most $K$ and identical rows within each community, whereas $P$ is the actual mean of the loop-free observations. Their difference is only the displayed diagonal correction. The matrix $W$ records the fluctuations of the observed graph about its mean.

Throughout the paper, generic constants may depend on $K,c_0,Q$ and a fixed clustering approximation factor, but not on $n$ or the label vector. They can be chosen uniformly when $Q$ ranges over a compact family with a common positive entry lower bound and a common positive minimum Euclidean row separation.

\subsection{Permutation-invariant losses}

Community names have no intrinsic meaning: exchanging the names of two communities leaves the partition unchanged. We therefore choose the best common relabeling before measuring error. For a soft assignment, the contribution of vertex $i$ is the probability mass placed outside its true community; for a hard assignment, it is the indicator of an incorrect label.
Let $\Simp=\{\tau\in[0,1]^{n\times K}:\tau\one=\one\}$. Define
\begin{align}
 e(\tau)&=\min_{\sigma\in\Perm}\frac1n\sum_i\{1-\tau_{i,\sigma(z_i^*)}\},
 &\cB(r)&=\{\tau:e(\tau)\le r\},\label{eq:softloss}\\
 R(\widehat z,z^*)&=\min_{\sigma\in\Perm}\frac1n\sum_i\one\{\widehat z_i\ne\sigma(z_i^*)\}.
 \label{eq:hardloss}
\end{align}
For an aligned pair of membership matrices, put
\begin{equation}\label{eq:metric}
 \ell(\tau,\tau')=\frac1{2n}\normone{\tau-\tau'},
\end{equation}
where $\normone{\cdot}$ denotes the entrywise sum of absolute values. Under correct alignment, $e(\tau)=\ell(\tau,Z^*)$. We use $\normop{\cdot}$, $\normF{\cdot}$ and $\normmax{\cdot}$ for operator, Frobenius and maximum-entry norms.

For example, a correctly assigned hard row contributes zero to $e$, a wrongly assigned hard row contributes $1/n$, and a row putting probability $0.8$ on its true class contributes $0.2/n$. Thus the same loss records both the initial hard clustering errors and the uncertainty in later soft iterates. The factor $1/(2n)$ in \eqref{eq:metric} gives exactly this interpretation because mass removed from the true coordinate reappears in the other coordinates.

\begin{proposition}[Soft-to-hard rounding]\label{prop:rounding}
If $\widehat z_i\in\argmax_k\tau_{ik}$, then
\begin{equation}\label{eq:rounding}
 R(\widehat z,z^*)\le2e(\tau).
\end{equation}
\end{proposition}

\begin{remark}[Why average loss is used]
If a hard initializer misclassifies even one vertex, its maximum-entry distance from the aligned truth is $1$. Its average error can nevertheless be small when few vertices are wrong. Average loss therefore matches the guarantee supplied by spectral initialization. It also permits exceptional vertices with large perturbations by controlling their contribution to the average loss in expectation, rather than imposing a uniform row-wise bound over all vertices.
\end{remark}

\subsection{Nodewise information}

To understand the statistical target, temporarily reveal all labels except that of one vertex. Conditional on these revealed labels, testing whether this vertex belongs to \(k\) or \(\ell\) amounts to comparing two product Bernoulli distributions for its incident edges. There are $n_m-\one\{m=k\}$ observed neighbors in community $m$ when its true label is $k$. The following definitions record the information in this idealized test; the implemented procedure never receives the revealed labels.
For $s\in[0,1]$, define
\begin{equation}\label{eq:hs}
 h_s(p,q)=-\log\{p^{1-s}q^s+(1-p)^{1-s}(1-q)^s\}.
\end{equation}
For $k\ne\ell$, let
\begin{align}
 J_{k\ell,n}&=\max_{0\le s\le1}\sum_{m=1}^K
 (n_m-\one\{m=k\})h_s(\rho_nq_{km},\rho_nq_{\ell m}),\label{eq:Jpair}\\
 J_n&=\min_{k\ne\ell}J_{k\ell,n}.\label{eq:J}
\end{align}
The same value of $s$ is used for all terms of a given pairwise test. For a fixed pair $k\ne\ell$, $J_{k\ell,n}$ is the Chernoff information for testing a vertex's label from its incident edges when all other labels are known. In the classification rule, the mixing proportions contribute the log odds $\log(\pi_{k,n}^*/\pi_{\ell,n}^*)$. Under Assumption~\ref{ass:proportions}, this term shifts the test threshold by only a bounded amount and hence does not change the first-order error exponent.

\begin{lemma}[Information expansion]\label{lem:information}
Define
\begin{equation}\label{eq:CH}
 \cC(\pi,Q)=\min_{k\ne\ell}\max_{0\le s\le1}
 \sum_m\pi_m\{(1-s)q_{km}+sq_{\ell m}-q_{km}^{1-s}q_{\ell m}^{s}\}.
\end{equation}
Under Assumptions~\ref{ass:proportions}--\ref{ass:sparsity}, uniformly in the labels,
\begin{equation}\label{eq:Jexpand}
 J_n=d_n\cC(\pistar,Q)+O(d_n\rho_n+\rho_n).
\end{equation}
There is a constant $\kappa>0$ such that
\begin{equation}\label{eq:Jbounds}
 \kappa d_n\le J_n\le2D_n.
\end{equation}
\end{lemma}

\begin{remark}[Unequal community sizes]
A vertex has approximately $n\pi_m$ potential neighbors in community $m$. The evidence that distinguishes labels $k$ and $\ell$ must therefore be weighted by every $\pi_m$, as in \eqref{eq:CH}, rather than just by the smallest proportion. For each fixed pair $k\ne\ell$, the maximized sum in \eqref{eq:CH} is the Chernoff--Hellinger divergence between the $K$-dimensional weighted connectivity profiles $(\pi_1q_{k1},\ldots,\pi_Kq_{kK})$ and $(\pi_1q_{\ell1},\ldots,\pi_Kq_{\ell K})$.
\end{remark}

\begin{remark}[Balanced planted partition as a check]
For $\pi_m=1/K$ and
$Q=(a-b)I_K+b\one\one^\top$, with $a>b>0$,
the general information quantity reduces to the classical planted-partition
Chernoff information. Specifically,
\[
\mathcal C(\pi,Q)=\frac{(\sqrt a-\sqrt b)^2}{K},
\]
and if
\[
I_n=-2\log\{\sqrt{a\rho_n b\rho_n}
+\sqrt{(1-a\rho_n)(1-b\rho_n)}\},
\]
then
\[
J_n=nI_n/K+O(\rho_n).
\]
Hence the normalization $nI_n/2$ corresponds specifically to the two-community
case.
\end{remark}

\subsection{Spectral initialization and simultaneous variational EM}

The two stages have different tasks. Spectral clustering supplies a partition with a small average error; VEM refines that partition while learning its nuisance parameters. Only the spectral stage uses the trimmed graph. Every subsequent likelihood score uses the original adjacency matrix, so trimming does not permanently discard a vertex's incident-edge information.
Let $d_i=\sum_jA_{ij}$ and $\bar d=n^{-1}\sum_i d_i$. Set
\begin{align}
 c_{\rm tr}&=200c_0^{-2},\qquad
 v_i=\one\{d_i\le c_{\rm tr}\max(\bar d,1)\},\notag\\
 A_{\rm tr}&=\diag(v)A\diag(v).\label{eq:trim}
\end{align}
Take the $K$ eigenvalues of $A_{\rm tr}$ with largest absolute values and their eigenvectors, denoted by $\widehat\Lambda_K$ and $\widehat U_K$. Apply an $\alpha$-approximate $K$-means procedure, with fixed $\alpha$, to the rows of
\begin{equation}\label{eq:embedding}
 X=\widehat U_K\widehat\Lambda_K.
\end{equation}
The approximation is relative to the optimal squared Euclidean clustering objective. This is the symmetric SC-RRE embedding of \citet[Algorithm~4]{ZhouAmini2019}, with the regularization threshold specified above. No row normalization is performed. If $\widehat z^{(0)}$ is the resulting assignment, initialize
\begin{equation}\label{eq:init}
 \tau^{(0)}_{ik}=(1-n^{-1})\one\{\widehat z_i^{(0)}=k\}+\frac1{nK}.
\end{equation}

Starting from $\tau^{(0)}$, the following steps define the simultaneous VEM
refinement on the original graph. Given the current soft assignment matrix
$\tau$, the M-step estimates the nuisance parameters through the following soft
counts.

The product $\tau_{ik}\tau_{j\ell}$ acts as a soft indicator that the ordered
pair $(i,j)$ belongs to block $(k,\ell)$, replacing the hard indicator
available under known labels. Accordingly, $N_{k\ell}$ is the resulting soft
edge count and $H_{k\ell}$ is the corresponding number of admissible ordered
dyads after excluding self-pairs. The subtraction in $H_{k\ell}$ removes the
inadmissible pairs with $i=j$.
\begin{align}
 s_k(\tau)&=\sum_i\tau_{ik},&
 \widehat\pi_k(\tau)&=s_k(\tau)/n,\label{eq:pimstep}\\
 N_{k\ell}(\tau)&=\sum_{i\ne j}\tau_{ik}\tau_{j\ell}A_{ij},&
 H_{k\ell}(\tau)&=s_ks_\ell-\sum_i\tau_{ik}\tau_{i\ell},\label{eq:softcounts}\\
 \widehat B_{k\ell}(\tau)&=\clip\!\left(\frac{N_{k\ell}(\tau)}{H_{k\ell}(\tau)},n^{-2},1-n^{-2}\right).
 &&\label{eq:Bmstep}
\end{align}
The symmetry of $A$ makes $\widehat B$ symmetric. Write $\phi(x,p)=x\log p+(1-x)\log(1-p)$.

The score for a candidate class adds its estimated log proportion to the expected edge log likelihood under the current neighbor assignments. Softmax converts these scores to probabilities, giving larger mass to classes that better explain the observed row of the graph.

The contrast between two scores, rather than either score alone, determines this choice: the ratio of their updated probabilities is $\exp\{L_{ik}(\tau)-L_{i\ell}(\tau)\}$. Consequently, the proof will control how estimation errors reduce the true class's score advantage. In a simultaneous update every row uses the previous membership matrix; newly computed rows are not substituted into later scores within the same iteration.
\begin{align}
 L_{ik}(\tau)&=\log\widehat\pi_k(\tau)+\sum_{j\ne i}\sum_m\tau_{jm}
 \phi(A_{ij},\widehat B_{km}(\tau)),\label{eq:score}\\
 \cT_A(\tau)_{ik}&=\frac{\exp L_{ik}(\tau)}{\sum_{\ell=1}^K\exp L_{i\ell}(\tau)}.
 \label{eq:update}
\end{align}
Although the local neighborhoods used in the results allow membership
matrices with zero entries, their closeness to the true membership matrix and
the lower bound on community proportions ensure that all class masses and
M-step denominators remain positive. The actual VEM iterates are strictly
positive in exact arithmetic because the softmax update assigns positive mass
to every class.

\begin{algorithm}
\caption{Spectral--VEM with unknown connectivity and proportions}
\label{alg:vem}
\begin{algorithmic}[1]
\Require $A$, fixed $K$, and a valid lower bound $c_0$.
\State Form $A_{\rm tr}$ by \eqref{eq:trim}; compute $X$ in \eqref{eq:embedding}.
\State Apply fixed-factor approximate $K$-means and initialize by \eqref{eq:init}.
\State Set $T_n=\max\{1,\lceil64c_0^{-2}\bar d/\log(\econst+\bar d)\rceil\}$.
\For{$t=0,\ldots,T_n-1$}
  \State Compute $\widehat\pi^{(t)}=\widehat\pi(\tau^{(t)})$ and
  $\widehat B^{(t)}=\widehat B(\tau^{(t)})$ using the original $A$.
  \State Compute all scores in \eqref{eq:score} using the same $\tau^{(t)}$.
  \State Set $\tau^{(t+1)}=\cT_A(\tau^{(t)})$ by simultaneous softmax.
\EndFor
\State Set $\widehat z_i\in\argmax_k\tau^{(T_n)}_{ik}$ and recompute
 $\widehat B=\widehat B(\tau^{(T_n)})$, $\widehat\pi=\widehat\pi(\tau^{(T_n)})$.
\Ensure $\widehat z$, $\tau^{(T_n)}$, $\widehat B$, and $\widehat\pi$.
\end{algorithmic}
\end{algorithm}

\begin{proposition}[Variational objective and update identities]\label{prop:objective}
For a fixed membership matrix \(\tau\) with positive class masses and \(H_{k\ell}(\tau)>0\) for every \(k,\ell\), the M-step \eqref{eq:pimstep}--\eqref{eq:Bmstep} maximizes
\begin{equation}\label{eq:ELBO}
\begin{split}
 \mathcal F(\tau,B,\pi)
 ={}&\sum_{i<j}\sum_{k,\ell}\tau_{ik}\tau_{j\ell}\phi(A_{ij},B_{k\ell})\\
 &+\sum_{i,k}\tau_{ik}\log\pi_k-\sum_{i,k}\tau_{ik}\log\tau_{ik},
\end{split}
\end{equation}
over symmetric $B$ with entries in $[n^{-2},1-n^{-2}]$ and over the probability simplex for $\pi$. For fixed \(B,\pi\) and the remaining membership rows, the softmax update \eqref{eq:update} is the unique maximizer over the current membership row.
\end{proposition}

\begin{remark}[The variational objective and the absence of count constraints]
The model conditions on the realized labels, but the algorithm does not
use or enforce the realized community counts. Thus \eqref{eq:ELBO} serves as
a working variational criterion for the deterministic-label model; if the
labels are instead generated independently from a categorical distribution,
it becomes the usual mean-field variational lower bound.
\end{remark}

\begin{remark}[Information used by algorithm]
The quantities \(q_-,q_+,\rho_n\), and \(J_n\) are used only for theoretical analysis and
are not supplied to the algorithm. The iteration budget depends only on the
observed average degree \(\bar d\) and a valid lower bound \(c_0\), not on the
unknown community proportions themselves. The constants in the theoretical
budget are chosen for the uniform asymptotic argument.
\end{remark}

\section{Statistical Guarantees for the VEM Iterates}
\label{sec:mainresults}
The next three theorems separate initialization, local stability, and expected statistical accuracy. Their assertions are uniform over deterministic labels satisfying Assumption~\ref{ass:proportions}. The proofs use one graph throughout. The principal conclusion is Theorem~\ref{thm:refinement}: after the prescribed iteration budget, the returned VEM iterate retains the oracle nodewise Chernoff exponent even though both nuisance parameters and soft labels are learned from that graph.

The order of the results is important. The first theorem places the initial error inside a fixed neighborhood. The second then supplies a recursion that is valid for every iterate in that neighborhood. The third controls the average over all graphs, including those on which initialization or local stability fails; it is not a direct consequence of discarding those graphs.

\begin{theorem}[Spectral initialization enters the basin]\label{thm:spectral}
Under Assumptions~\ref{ass:proportions}--\ref{ass:sparsity}, there are constants $C>0$ and an event $\Omega_n$ with
\begin{equation}\label{eq:OmegaProb}
 \Prob(\Omega_n^c)\le Ce^{-10D_n}
\end{equation}
such that the initializer in Algorithm~\ref{alg:vem} satisfies, on $\Omega_n$,
\begin{align}
 R(\widehat z^{(0)},z^*)&\le C/d_n,\label{eq:spectralrate}\\
 e(\tau^{(0)})&=(1-n^{-1})R(\widehat z^{(0)},z^*)+n^{-1}(1-1/K)
 \le C/d_n+1/n.\label{eq:initialsoft}
\end{align}
In particular, it enters every fixed positive average-loss neighborhood for sufficiently large $n$.
\end{theorem}
\begin{remark}[Rank deficiency and negative signal eigenvalues]
Distinct rows of $Q$ give distinct community row centers even when $Q$ has deficient rank. This is the reduced-rank geometry of \citet{ZhouAmini2019}. Weighting eigenvectors by their eigenvalues keeps small extra directions from having the same influence as signal directions; ordering by magnitude retains negative signal eigenvalues. The exponentially small initialization failure bound is needed because failure also contributes to the expected risk in Theorem~\ref{thm:refinement}.
\end{remark}

\begin{theorem}[Uniform local loss contraction]\label{thm:contraction}
Under Assumptions~\ref{ass:proportions}--\ref{ass:sparsity}, there are constants $r_0,C,c>0$ such that, with probability at least $1-Ce^{-cd_n^{3/4}}$,
\begin{equation}\label{eq:mainrecursion}
 e\{\cT_A(\tau)\}\le\gamma_ne(\tau)+\varepsilon_n
 \quad\text{simultaneously for every }\tau\in\cB(r_0),
\end{equation}
where
\begin{equation}\label{eq:gammaepsilon}
 \gamma_n=Cd_n^{-1/2}<1,\qquad
 \varepsilon_n=C\exp\{-J_n+Cd_n^{3/4}\}.
\end{equation}
On this event, every trajectory started with $e_0\le r_0$ remains in the basin and obeys
\begin{equation}\label{eq:iterated}
 e_t\le\gamma_n^te_0+\frac{1-\gamma_n^t}{1-\gamma_n}\varepsilon_n,
 \qquad t\ge0.
\end{equation}
\end{theorem}
\begin{remark}[Uniformity is the same-graph guarantee]
Every iterate is computed from $A$, so one cannot apply a concentration bound that assumes its input is independent of $A$. Equation~\eqref{eq:mainrecursion} instead holds for every admissible input on one graph event. It therefore applies successively to the random iterates, without an additional failure probability for each iteration. This controls their errors relative to the truth; it does not assert that the empirical map contracts the distance between two arbitrary inputs.
\end{remark}

\begin{theorem}[Expected statistical refinement]\label{thm:refinement}
Under Assumptions~\ref{ass:proportions}--\ref{ass:sparsity}, Algorithm~\ref{alg:vem}, with
\begin{equation}\label{eq:iterations}
 T_n=\max\left\{1,\left\lceil
 \frac{64c_0^{-2}\bar d}{\log(\econst+\bar d)}\right\rceil\right\},
\end{equation}
satisfies
\begin{align}
 \E e(\tau^{(T_n)})&\le C\exp\{-J_n+Cd_n^{3/4}\},\label{eq:expectedsoft}\\
 \E R(\widehat z,z^*)&\le C\exp\{-J_n+Cd_n^{3/4}\}
 =\exp\{-(1-o(1))J_n\}.\label{eq:expectedhard}
\end{align}
On an event of probability tending to one, $T_n=O(d_n/\log d_n)$.
\end{theorem}

\begin{remark}[High-probability contraction does not directly imply the expected-risk exponent]
The contraction bound in Theorem~\ref{thm:contraction} holds with failure
probability of order $e^{-cd_n^{3/4}}$. Simply averaging
\eqref{eq:iterated} and using the trivial bound $e_t\le1$ outside that event
would therefore leave a contribution of this order, which is larger than the
target $e^{-J_n+o(d_n)}$. The proof of Theorem~\ref{thm:refinement} instead
retains the common random remainder from Lemma~\ref{lem:onestep} and controls
its expectation directly. This expectation bound simultaneously controls the
accumulated error inside the shrinking neighborhood and the probability of
leaving that neighborhood, thereby preserving the desired expected-risk
exponent.
\end{remark}

\begin{corollary}[Vanishing error recovery and nuisance parameter consistency]\label{cor:partial}
Under the preceding assumptions,
\begin{equation}\label{eq:partial}
 R(\widehat z,z^*)\longrightarrow0\quad\text{in probability}.
\end{equation}
After a common label alignment,
\begin{align}
 \normmax{\widehat\pi-\pistar}&\le e(\tau^{(T_n)}),\label{eq:pierror}\\
 \frac{\normmax{\widehat B-B^*}}{\rho_n}&=O_{\Prob}(d_n^{-1/2}).\label{eq:Berror}
\end{align}
\end{corollary}
\begin{remark}[Recovery regimes and nuisance-parameter estimation]
The convergence $R(\widehat z,z^*)\to0$ means that the fraction of
misclassified vertices vanishes, although the number of errors may remain
positive; this is often called almost exact recovery and does not by itself
imply exact recovery. The proportion estimate in \eqref{eq:pierror} targets
the realized class fractions under the fixed-label formulation, which differ
from population proportions under a random-label model due to sampling
fluctuations. The rate in \eqref{eq:Berror} is sufficient for the variational
refinement analysis and is not intended as a sharp parametric rate for
estimating $B^*$.
\end{remark}

The preceding results establish an upper bound for the estimator returned by
Algorithm~1. To establish optimality of this recovery exponent, one also
needs a matching lower bound over a precisely defined parameter space. The
next section compares the achieved misclassification risk with the minimax
risk and then derives the sharp threshold for exact recovery.

\section{Statistical Optimality of the Recovery Procedure}\label{sec:optimality}

To establish statistical optimality, we compare the risk of the estimator
returned by Algorithm~\ref{alg:vem} with the minimax risk over a local
parameter space that allows both the connectivity matrix and community
proportions to vary.

Fix a symmetric positive matrix $Q^0$ with distinct rows and a probability vector $\pi^0$ with $\min_k\pi_k^0>c_0$. Let
\begin{equation}\label{eq:delta}
 \delta_n\to0,\qquad \frac{n\delta_n^2}{d_n}\to\infty.
\end{equation}
For example, $\delta_n=\rho_n^{1/4}$ is admissible. Define
\begin{equation}\label{eq:Theta}
 \Theta_n=\{(Q,z):Q=Q^\top,\ \normmax{Q-Q^0}\le\delta_n,
 \ \normmax{\pi(z)-\pi^0}\le\delta_n\},
\end{equation}
where the graph has connection matrix $\rho_nQ$, and put
\begin{equation}\label{eq:Gamma}
 \Gamma_n=d_n\cC(\pi^0,Q^0).
\end{equation}
For sufficiently large $n$, every element of this class satisfies the
positivity and separation assumptions uniformly. Neither the realized
$Q$ nor the realized community counts are supplied to an estimator.

The neighborhood shrinks around a fixed profile, but it still contains many different label configurations. Condition~\eqref{eq:delta} ensures that the neighborhood is large enough
to accommodate the local alternatives used in the lower-bound construction,
while still shrinking to the reference profile. The reference points $Q^0$ and $\pi^0$ define the local parameter neighborhood over which the minimax comparison is made; they are not inputs
to Algorithm~\ref{alg:vem}.

Allowing the counts to vary is part of the statistical comparison.
It permits nearby models to differ in one vertex's label. If exact class
sizes were imposed instead, such an alternative might be inadmissible and
a different lower-bound argument would be needed. The minimax theorem below
uses precisely the space in \eqref{eq:Theta}.

\begin{theorem}[Minimax misclassification exponent]\label{thm:minimax}
Under \eqref{eq:delta} and Assumption~\ref{ass:sparsity},
\begin{equation}\label{eq:minimax}
 \inf_{\widehat z}\sup_{(Q,z)\in\Theta_n}
 \E_{Q,z}R(\widehat z,z)
 =\exp\{-(1+o(1))\Gamma_n\}.
\end{equation}
\end{theorem}

The estimator returned by Algorithm~\ref{alg:vem} attains this minimax
exponent over $\Theta_n$. Writing $\mathfrak R_n^*$ for the minimax risk in \eqref{eq:minimax}, the assertion is equivalent to
\[
 \frac{\log\mathfrak R_n^*}{\Gamma_n}\longrightarrow-1.
\]
Thus the minimax optimality is first-order on the logarithmic scale.

\begin{remark}[Why the count constraint matters]
Changing one vertex from class $k$ to class $\ell$ changes two class counts by one.
Such local alternatives are allowed in \eqref{eq:Theta}, where counts may vary,
but may be inadmissible if exact counts are prescribed. The nodewise lower bound
uses this freedom, so the minimax conclusion here does not apply to a parameter
space with known exact counts.
\end{remark}

Expected loss measures the average fraction of misclassified vertices,
whereas exact recovery requires that no vertex be misclassified.
The bound $\one\{R>0\}\le nR$ shows that a sufficiently small expected
loss implies exact recovery. Conversely, a lower bound on expected loss
alone does not imply that an error occurs with high probability, so the
exact-recovery converse requires a separate argument.

\begin{theorem}[Sharp first-order exact recovery threshold]\label{thm:exact}
For deterministic labels satisfying Assumption~\ref{ass:proportions}, if, for fixed $\epsilon>0$,
\begin{equation}\label{eq:above}
 J_n\ge(1+\epsilon)\log n,
\end{equation}
then Algorithm~\ref{alg:vem} satisfies $\Prob\{R(\widehat z,z^*)=0\}\to1$.
For the class \eqref{eq:Theta}, if
\begin{equation}\label{eq:below}
 \Gamma_n\le(1-\epsilon)\log n,
\end{equation}
then
\begin{equation}\label{eq:strongconverse}
 \inf_{\widehat z}\sup_{(Q,z)\in\Theta_n}
 \Prob_{Q,z}\{R(\widehat z,z)>0\}\longrightarrow1.
\end{equation}
\end{theorem}

The converse is proved by showing that below the threshold, with high probability at least one vertex admits a more likely one-label alternative

\begin{corollary}[Independent categorical labels]\label{cor:random}
Suppose instead that $z_i^*$ are independent with distribution $\Cat(\pi)$, where $\min_k\pi_k>c_0$, and retain the connectivity and sparsity assumptions. Then
\begin{equation}\label{eq:randomrisk}
 \E R(\widehat z,z^*)\le\exp\{-(1-o(1))d_n\cC(\pi,Q)\}.
\end{equation}
Moreover, with label alignment,
\[
 \normmax{\widehat\pi-\pi}
 \le e(\tau^{(T_n)})+\normmax{\pi_n^*-\pi},
 \qquad \normmax{\pi_n^*-\pi}=O_{\Prob}(n^{-1/2}).
\]
Exact recovery holds when
$d_n\cC(\pi,Q)\ge(1+\epsilon)\log n$.
If
$d_n\cC(\pi,Q)\le(1-\epsilon)\log n$,
then even an estimator supplied with the true parameters has
exact recovery probability tending to zero under this label law.
\end{corollary}
\begin{remark}
Conditional on the realized labels, the deterministic-label result applies
whenever the empirical community proportions remain bounded away from zero.
To average the exponential risk bound over the random labels, we also need
the probability that the realized proportions deviate from $\pi$ to be
negligible on the same risk scale. A suitably shrinking neighborhood around
$\pi$ provides this control.
\end{remark}

\section{Mechanism of the Proof}\label{sec:mechanism}
The central estimate is Lemma~\ref{lem:onestep}, which gives a uniform
one-step bound with a common random remainder whose expectation is controlled
at the Chernoff scale. This estimate underlies both the high-probability local
contraction and the unconditional expected-risk result. The population fixed-point result describes the local geometry, but it is not used in the proof of the empirical expected-risk bound. Complete proofs and supporting details are given in Sections~C--G of the supplementary material.

The proof proceeds in three steps. First, we construct a single high-probability
graph event on which the bilinear and operator-norm bounds needed below hold
simultaneously for all admissible soft membership matrices. These uniform
bounds then control the block-probability and mixing-proportion estimates at
every data-dependent iterate. Finally, for each vertex and each competing
community, we compare the actual VEM score difference with the oracle
likelihood-ratio score obtained from the true parameters and labels.

This comparison does not require every vertex--class score perturbation to be
small. Instead, perturbations exceeding a prescribed margin are controlled in
aggregate. For vertices retained after trimming, the aggregate contribution of
these large perturbations is bounded by a multiple of the current membership
error, yielding the contraction term. Cases where the oracle score itself has a
small margin, together with the contribution of unusually high-degree vertices,
are absorbed into a common random remainder. Controlling the expectation of
this remainder, rather than replacing it immediately by a failure event, is what
preserves the Chernoff exponent in the final expected-risk bound.

\begin{lemma}[A uniform graph event]\label{lem:graph}
There is an event $\Omega_n$ satisfying \eqref{eq:OmegaProb} on which
\begin{align}
 c_0^2D_n/2&\le\bar d\le2D_n,\label{eq:average}\\
 \sup_{x,y\in[0,1]^n}|x^\top Wy|&\le Cn\sqrt{d_n},\label{eq:bilinear}\\
 \normop{A_{\rm tr}-P}&\le C\sqrt{d_n}.\label{eq:operator}
\end{align}
If
\begin{equation}\label{eq:V}
 V_n=\sum_i(d_i+D_n)\one\{d_i>100D_n\},
\end{equation}
then
\begin{equation}\label{eq:Vmean}
 \E V_n\le CnD_ne^{-90D_n}.
\end{equation}
\end{lemma}
\begin{remark}[Roles of the three graph controls]
The M-step numerator contains
$\tau_{\cdot k}^{\top}W\tau_{\cdot\ell}$, so the bilinear bound must
hold uniformly over all soft weight vectors, including data-dependent ones.
In the E-step, the operator-norm bound controls the noise--membership
interaction on the vertices retained after trimming through its aggregate
squared size. The contribution of the removed high-degree vertices is handled
separately through $V_n$ and controlled in expectation. The three controls
therefore address different sources of error in the update analysis.
\end{remark}

On the graph event, soft class sizes cannot be too small when the average membership error is small. Dividing soft edge counts by these well-behaved denominators therefore gives the next estimate. It separates error due to the current labels from error due to the random graph.

\begin{lemma}[Uniform M-step control]\label{lem:Mstep}
There is a fixed $r_0>0$ such that, on $\Omega_n$, for every correctly aligned $\tau$ with $e(\tau)\le r_0$,
\begin{align}
 \normmax{\widehat\pi(\tau)-\pistar}&\le e(\tau),\label{eq:uniformpi}\\
 \frac{\normmax{\widehat B(\tau)-B^*}}{\rho_n}
 &\le C\{e(\tau)+d_n^{-1/2}\}.\label{eq:uniformB}
\end{align}
The denominators in \eqref{eq:Bmstep} are at least $cn^2$, and the clipping is inactive for sufficiently large $n$.
If $\widehat B_P$ denotes the same update with $P$ in place of $A$, then in the aligned neighborhood
\begin{align}
 \normmax{\widehat B_P(\tau)-\widehat B_P(\tau')}&\le C\rho_n\ell(\tau,\tau'),\label{eq:BPlip}\\
 \normmax{\widehat\pi(\tau)-\widehat\pi(\tau')}&\le\ell(\tau,\tau').\label{eq:pilip}
\end{align}
\end{lemma}
\begin{remark}[Effect of estimating mixing proportions]
The proportion estimate differs from the truth by at most $e(\tau)$, and all
relevant proportions remain bounded away from zero. The resulting perturbation
of a mixing log-ratio is therefore only $O(e(\tau))$. By contrast, errors in
estimating $B$ accumulate over the neighbors of a vertex, producing a
perturbation of order $O(d_ne(\tau)+\sqrt{d_n})$ on retained vertices.
The contribution from estimating the mixing proportions is thus absorbed into
this larger score perturbation and does not affect the leading information
exponent.
\end{remark}

Because community labels are only identifiable up to permutation, we fix one
labeling under which the estimated and true communities are aligned, and define
the closed convex neighborhood
\[
 \mathcal B_{\mathrm{id}}(r_0)
 =\{\tau\in\Simp:\ell(\tau,Z^*)\le r_0\}.
\]
The uniqueness statement below is with respect to this fixed alignment.
By equivariance, the same conclusion holds for the corresponding neighborhoods
under other label permutations.

\begin{proposition}[Local contraction of the population map]
\label{prop:population}
Let $\cT_P$ be the map obtained by replacing $A$ by $P$ in both the
M-step and the E-step. For $r_0>0$ sufficiently small but fixed,
$\cT_P$ maps $\mathcal B_{\mathrm{id}}(r_0)$ into itself and satisfies
\begin{equation}\label{eq:popcontraction}
 \ell\{\cT_P(\tau),\cT_P(\tau')\}
 \le\beta_n\ell(\tau,\tau'),\qquad
 \beta_n=Cd_ne^{-cd_n}<1,
\end{equation}
for all $\tau,\tau'\in\mathcal B_{\mathrm{id}}(r_0)$.
Consequently, $\cT_P$ has a unique fixed point
$\tau_P^\dagger\in\mathcal B_{\mathrm{id}}(r_0)$, and
\begin{equation}\label{eq:popbias}
 \ell(\tau_P^\dagger,Z^*)
 \le\frac{\ell\{\cT_P(Z^*),Z^*\}}{1-\beta_n}
 \le Ce^{-cd_n}.
\end{equation}
\end{proposition}

\begin{remark}[The hard truth is not a soft fixed point]
Finite scores give strictly positive softmax probabilities in every class,
whereas $Z^*$ assigns zero probability to all wrong classes. Hence
$\cT_P(Z^*)\ne Z^*$ at finite $n$, so a zero-residual contraction around
$Z^*$ cannot hold. Proposition~\ref{prop:population} instead identifies a
nearby soft fixed point. This population result is used only to describe the
local geometry. For the observed graph, the one-step bound is proved directly
by comparing empirical score differences with oracle likelihood-ratio scores
and controlling the resulting perturbations in aggregate.
\end{remark}

We now return to the observed graph. Rather than subtracting a population
softmax from an empirical softmax, we compare their underlying classification
evidence with an ideal score computed at the truth. Positive oracle contrasts
favor the true class. The margin $h$ below allows for the discrepancy between
the fitted VEM score and the oracle score caused by parameter estimation and
imperfect neighbor memberships.

For $z_i^*=k$ and $\ell\ne k$, define the oracle score
\begin{equation}\label{eq:oracle}
 G_{i\ell}=\log\frac{\pi_{k,n}^*}{\pi_{\ell,n}^*}
 +\sum_{j\ne i}\{\phi(A_{ij},B^*_{k,z_j^*})-\phi(A_{ij},B^*_{\ell,z_j^*})\}.
\end{equation}
For $h>0$, put
\begin{align}
 F_h&=\frac1n\sum_i\sum_{\ell\ne z_i^*}\min\{1,e^{-G_{i\ell}+h}\},\label{eq:Fh}\\
 U_n(h)&=F_h+C\frac{V_n}{nh}.\label{eq:Uh}
\end{align}
Here $F_h$ measures the total wrong-class mass suggested by the oracle
comparisons after allowing each oracle score advantage to decrease by $h$.
The term $CV_n/(nh)$ accounts for the contribution of removed high-degree
vertices. A larger $h$ allows larger discrepancies between the fitted and
oracle scores, but weakens the exponential bound on the remainder. The next
lemma makes this trade-off explicit.

There are two sources of error after an update. First, the oracle score itself
may provide only weak separation between the true and a competing class; the
aggregate contribution of such cases is measured by $F_h$. Second, the fitted
scores may deviate substantially from the oracle scores. For vertices retained
after trimming, these large deviations are controlled in aggregate through a
squared-error bound, yielding a term proportional to the current membership
error. Contributions from removed high-degree vertices are handled separately
through $V_n$. This decomposition explains why the next bound contains both a
contraction term and a common additive remainder.

\begin{lemma}[Uniform empirical one-step bound]\label{lem:onestep}
On $\Omega_n$, for every aligned $\tau$ with $e(\tau)\le r_0$ and every $h$ satisfying
\begin{equation}\label{eq:hcondition}
 C_0\{d_ne(\tau)+\sqrt{d_n}\}\le h\le D_n,
\end{equation}
we have
\begin{equation}\label{eq:onestep}
 e\{\cT_A(\tau)\}\le C\frac{d_n}{h^2}e(\tau)+U_n(h).
\end{equation}
For every deterministic $h\ge C_0\sqrt{d_n}$,
\begin{equation}\label{eq:Umean}
 \E U_n(h)\le Ce^{-J_n+h}.
\end{equation}
\end{lemma}
\begin{remark}[Two roles of the margin]
The margin $h$ must cover the score perturbation
$C\{d_ne(\tau)+\sqrt{d_n}\}$, but increasing it weakens the residual bound
$e^{-J_n+h}$. On a fixed neighborhood, taking $h$ proportional to $d_n$
gives a coarse contraction that drives the error toward a smaller neighborhood.
Once $e(\tau)\le d_n^{-1/4}$, taking $h$ proportional to $d_n^{3/4}$ still
covers the perturbation, gives $d_n/h^2=O(d_n^{-1/2})$, and incurs only an
$o(d_n)$ loss in the exponent.
\end{remark}

\section{Computational Complexity}\label{sec:computation}

The VEM score includes contributions from both edges and nonedges, but
evaluating the nonedge terms does not require iterating over all $n^2$
dyads. The identity below first collects the contribution that would arise
if every potential dyad were a nonedge, and then adds a correction for the
dyads where an edge is observed. This gives an exact sparse implementation
of the theoretical score.

\begin{proposition}[Computational complexity and sparse implementation]\label{prop:complexity}
For fixed $K$ and $c_0$, Algorithm~\ref{alg:vem} can be implemented using
a polynomial number of operations. Given a sparse adjacency
representation, one M-step and one E-step together require
\begin{equation}\label{eq:complexity}
 O(|E|K+nK^2)
\end{equation}
arithmetic operations, apart from the initialization. With elementwise logarithms, let
\[
 G=\log(1-\widehat B),\qquad H_B=\log\{\widehat B/(1-\widehat B)\}.
\]
The score matrix has the exact identity
\begin{equation}\label{eq:sparseidentity}
 L=\one(\log\widehat\pi)^\top
 +(A\tau)H_B^\top+(\one s^\top-\tau)G^\top.
\end{equation}
\end{proposition}
\begin{remark}
For fixed $K$, exhaustive enumeration of centers chosen from the data points
requires polynomially many operations and guarantees the approximation factor
required by the initialization theorem. This shows that the statistical
guarantee can be attained by a procedure using polynomially many basic
operations.
\end{remark}

\section{Numerical Experiments}\label{sec:numerics}
We examine the practical behavior of the VEM iteration at finite sample sizes,
with unknown parameters and the same graph used at every update. The
experiments evaluate how VEM improves spectral initialization, how its
performance changes under different community proportions and connectivity
structures, whether it exhibits the predicted exact recovery transition, how
sensitive it is to initialization, and how the computational cost scales with
network size. The experiments also record settings in which the practical
implementation fails to show stable refinement. Each graph generation setting
is replicated 100 times.

\subsection{Simulation settings and information matching}

We set $K=3$ and use the following raw connectivity matrices:
\begin{align*}
 Q_{\rm raw}^{\rm A}&=\begin{pmatrix}6&1&1\\1&6&1\\1&1&6\end{pmatrix},&
 Q_{\rm raw}^{\rm H}&=\begin{pmatrix}7&0.5&1.5\\0.5&6.5&2\\1.5&2&5.5\end{pmatrix},\\
 Q_{\rm raw}^{\rm D}&=\begin{pmatrix}1&4&5\\4&2&4\\5&4&1\end{pmatrix},&
 Q_{\rm raw}^{\rm R}&=\begin{pmatrix}1.95&1&0.05\\1&1&1\\0.05&1&1.95\end{pmatrix}.
\end{align*}
Models A and H are homogeneous and heterogeneous assortative examples, respectively. Model D has two negative eigenvalues. Model R has rank two, while retaining positive entries and distinct rows. For $\pi^0=(0.5,0.3,0.2)$, we normalize each matrix once by
\[
 Q=Q_{\rm raw}/\{(\pi^0)^\top Q_{\rm raw}\pi^0\}.
\]
The normalized matrix is kept fixed when community proportions are varied. Except in the imbalance experiment, proportions equal $\pi^0$; in that
experiment, we consider the three proportion vectors
$(1/3,1/3,1/3)$, $(0.5,0.3,0.2)$, and $(0.7,0.2,0.1)$. Counts are rounded by the largest-remainder rule, and labels are randomly permuted over the vertices. Sampling a binomial number of dyads within each block and then a uniform subset of that size generates exactly the independent Bernoulli model.

For a target information value, $\rho_n$ is obtained by numerically solving
the finite-sample definition of $J_n$, using the common Chernoff parameter
for each ordered community pair. Thus different structures are compared at
matched information, not at matched raw edge probability. The resulting
connection probabilities remain in the sparse range, with
$\max_{k\ell}B^*_{k\ell}\le0.5$ in all settings.

\begin{table}[tbp]
\centering\small\setlength{\tabcolsep}{4pt}
\begin{tabular}{p{.15\textwidth}p{.53\textwidth}rr}\toprule
Experiment & Settings & Replicates & Graphs\\\midrule
Risk & A, H, D, R; $n=600,1800$; $J_n=1.5,3,4.5,6$ & 100 & 3,200\\
Exact recovery & A, H, D, R; $n=600,1800$; $J_n/\log n=0.6,0.8,1,1.2,1.4$ & 100 & 4,000\\
Imbalance & H; $n=1200$; $J_n=2,4,6$; three proportion vectors & 100 & 900\\
Initialization sensitivity & H, D, R; $n=1200$; $J_n=3,1.2\log n$; seven initial errors on each graph & 100 & 600\\
Timing & H; $n=600,1200,2400,4800$; $\rho_n=8\log(n)/n$ & 100 & 400\\\midrule
Total & 91 graph-generating settings & & 9,100\\\bottomrule
\end{tabular}

\vspace{2mm}
\caption{Replication counts refer to independent graphs per
graph-generating setting. Initial-error levels and comparisons among algorithmic
variants reuse the same graph realizations rather than increase the number of
independent networks.}\label{tab:design}
\end{table}

The initialization sensitivity experiment constructs controlled initial
conditions by corrupting a fraction
$\alpha\in\{0,0.05,0.1,0.2,0.3,0.4,0.5\}$
of each true community through cyclic relabeling, followed by the same
smoothing used for the spectral initialization. The truth is used only to
construct this diagnostic initialization and is not available to the estimator.

\subsection{Implementation, comparisons, and Monte Carlo uncertainty}

The practical estimator uses $c_0=0.08$, largest magnitude eigenpairs with
residual tolerance $10^{-9}$, and the weighted embedding without row
normalization. We use 30 $k$-means++ initializations followed by Lloyd
iterations and retain the solution with the smallest clustering objective.
The clustering step follows a standard practical implementation; the fixed
factor approximation condition in the initialization theorem is a theoretical
assumption rather than an empirically verified condition here. Every VEM step
estimates the full symmetric matrix $B$ and the proportion vector $\pi$ from
the current soft memberships, and updates all rows simultaneously on the observed
graph.

A run stops after three successive iterations, following at least ten
iterations, satisfy both a maximum membership change below $10^{-9}$ and a
parameter change criterion below $10^{-7}$. The parameter criterion is the
maximum of the relative Frobenius change in $B$ and the maximum absolute
change in $\pi$. We impose a cap of 200 iterations. Twenty additional updates
are then performed as a numerical check and do not replace the reported output.
This stopping rule differs from the conservative theoretical budget in
Algorithm~\ref{alg:vem}. 

The primary comparisons are spectral clustering alone (SC), one soft count
M/E update (One step), the practical Spectral VEM described above (SVEM), and
hard refinement using the same parameter estimates and score with an argmax
membership update (Hard). Hard is a matched variant of SVEM related to the
posterior thresholding update of \citet{LiLe2023}, rather than a reproduction
of their full procedure. Similarly, One step isolates the effect of a single
refinement update and is distinct from the blockwise procedure of
\citet{ZhouLi2018}.

On all risk replicates for H, D, and R at $n=1800$, we additionally compare
a known $B,\pi$ oracle procedure, five random start VEM runs selected by their
final working objective, and VEM initialized by the unweighted eigenvectors, to assess the contribution
of the proposed weighted spectral initialization.
For D we also use the $K$ largest algebraic, rather than largest magnitude,
eigenvalues. Oracle results are aligned by maximizing the known model working
objective over label permutations, not by using the true labels. The oracle
procedure serves only as a reference and is not an information theoretic lower
bound. Five start random VEM uses more computation than a single spectral
start.

We report graph-level means and standard errors estimated as $s/\sqrt{M}$,
where $s$ is the sample standard deviation and $M$ is the number of
independent graph realizations. Confidence intervals for mean performance are
obtained by percentile bootstrap with 2,000 resamples of the graph
realizations. Comparisons between methods use the same graph realizations and
are therefore paired. Exact recovery probabilities are reported with
two-sided exact 95\% binomial confidence intervals. Additional Monte Carlo
precision calculations and implementation details are provided in Section~K
of the supplementary material.

All numerical results were obtained on a MacBook Pro with an Apple M4
processor and 32 GB of memory, running macOS 26.6.2 and Python 3.12.13,
with NumPy 2.3.5, SciPy 1.17.0, and scikit-learn 1.8.0
\citep{HarrisEtAl2020,VirtanenEtAl2020,PedregosaEtAl2011}.
Timing measurements for all 100 graphs at each network size used one process
and one numerical thread with the Apple Accelerate numerical library.

\subsection{Refinement accuracy and exact recovery}

\begin{figure}[tbp]\centering
\includegraphics[width=.485\textwidth]{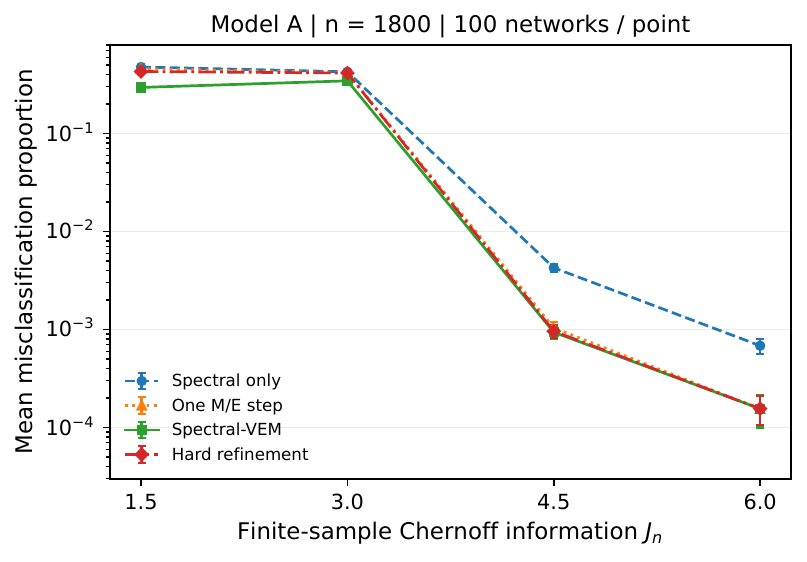}\hfill
\includegraphics[width=.485\textwidth]{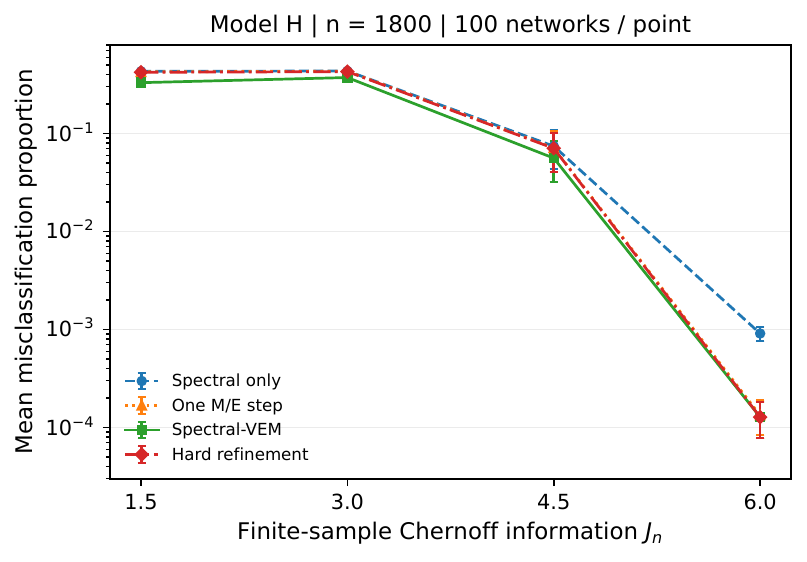}\\[3pt]
\includegraphics[width=.485\textwidth]{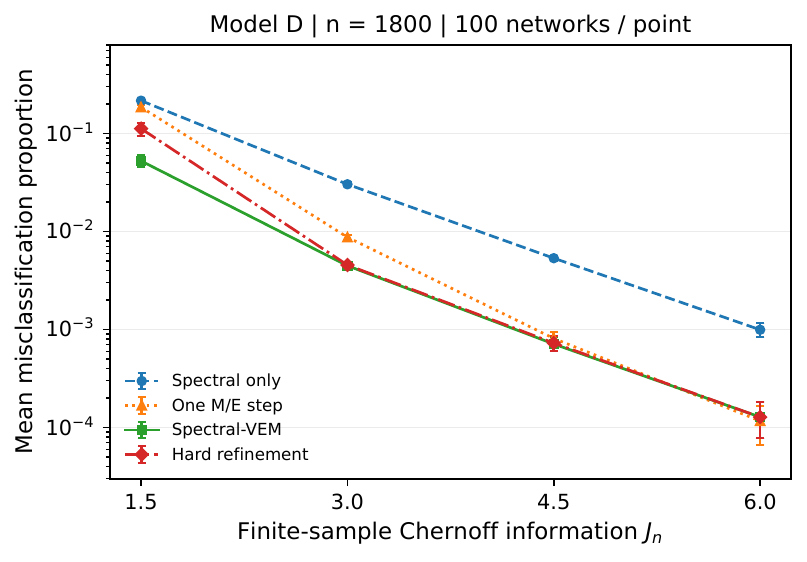}\hfill
\includegraphics[width=.485\textwidth]{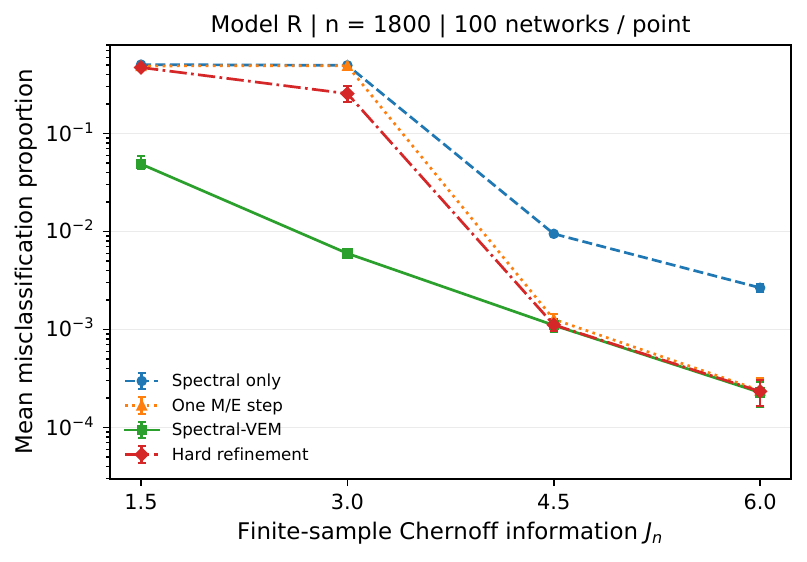}
\caption{Mean misclassification proportion at $n=1800$ and
$\pi^0=(0.5,0.3,0.2)$, with 95\% graph-level percentile bootstrap intervals.
Each point is based on 100 independent graphs, and all methods within a panel
are evaluated on the same graphs.}
\label{fig:risk}
\end{figure}

Figure~\ref{fig:risk} and Table~\ref{tab:highinfo} show clear reductions in
misclassification error relative to the spectral initializer at the higher
information levels. At $J_n=6$ and $n=1800$, SVEM reduces the mean error by
factors of approximately $4.4$, $7.1$, $7.8$, and $11.7$ for A, H, D, and R,
respectively. Table~\ref{tab:highinfo} gives the corresponding error levels
and the estimation errors for $B$ and $\pi$. One-step, SVEM, and hard
refinement give similar errors in several of these settings.

\begin{table}[tbp]
\centering\small
\setlength{\tabcolsep}{4pt}
\begin{tabular}{lrrrrrr}
\toprule
& \multicolumn{4}{c}{Mean $R\times10^3$ (Monte Carlo SE $\times10^3$)}
& \multicolumn{2}{c}{SVEM parameters}\\
\cmidrule(lr){2-5}\cmidrule(lr){6-7}
Model & SC & One-step & SVEM & Hard
& $\|\widehat B-B^*\|_F/\|B^*\|_F$
& $10^3\|\widehat\pi-\pi^*\|_\infty$\\
\midrule
A & 0.683 (0.062) & 0.156 (0.029) & 0.156 (0.030) & 0.156 (0.029) & 0.0152 & 0.155 \\
H & 0.911 (0.079) & 0.133 (0.027) & 0.128 (0.027) & 0.128 (0.027) & 0.0110 & 0.137 \\
D & 0.994 (0.082) & 0.117 (0.027) & 0.128 (0.027) & 0.128 (0.027) & 0.0068 & 0.125 \\
R & 2.667 (0.130) & 0.239 (0.038) & 0.228 (0.035) & 0.233 (0.036) & 0.0102 & 0.232 \\
\bottomrule
\end{tabular}
\vspace{1mm}
\caption{Recovery and nuisance parameter estimation at $n=1800$, $J_n=6$,
and $\pi^0=(0.5,0.3,0.2)$, using 100 independent graphs for each model.
Parameter errors are averaged after joint row and column label alignment.
Parentheses report Monte Carlo standard errors across graphs.}
\label{tab:highinfo}
\end{table}

At lower information levels, performance varies substantially across
connectivity structures. For H, the mean SVEM error at $J_n=3$ is $0.3725$,
and none of the 100 runs satisfies the stopping tolerances; at $J_n=4.5$ the
mean error is $0.0561$, with 84 of 100 runs satisfying the tolerances. By
contrast, D and R have mean errors $0.00448$ and $0.00598$ at $J_n=3$.
Thus, at these sample sizes, low-information refinement is sensitive to the
connectivity structure and can fail when the spectral initializer is not
sufficiently informative. Results for all 32 risk settings, including
$n=600$, are reported in Section~K of the supplementary material.

We next consider exact recovery, recording whether a graph is recovered with
no misclassified vertices. Since this is a binary outcome at the graph level,
uncertainty is summarized by exact binomial confidence intervals.

\begin{figure}[tbp]\centering
\includegraphics[width=.485\textwidth]{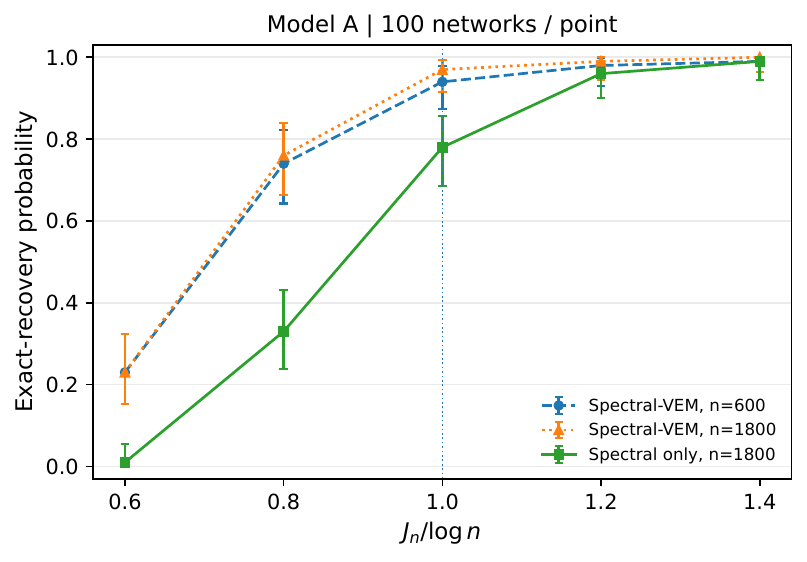}\hfill
\includegraphics[width=.485\textwidth]{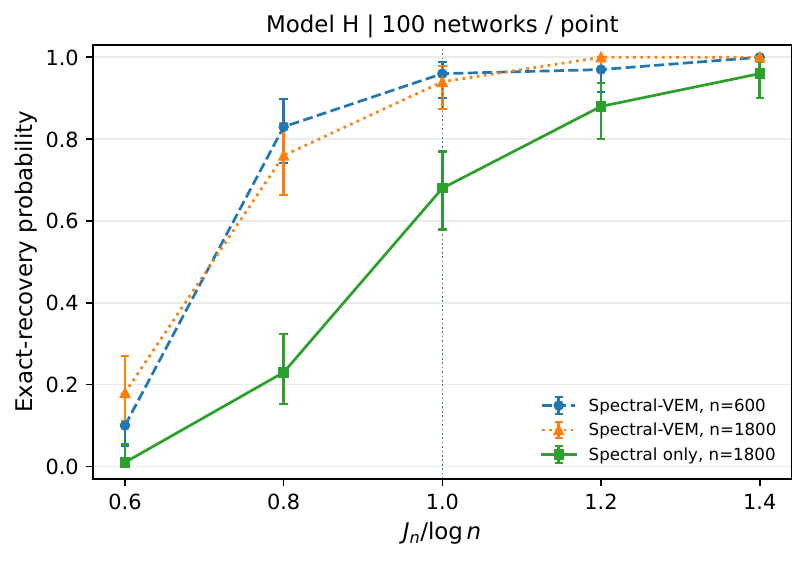}\\[3pt]
\includegraphics[width=.485\textwidth]{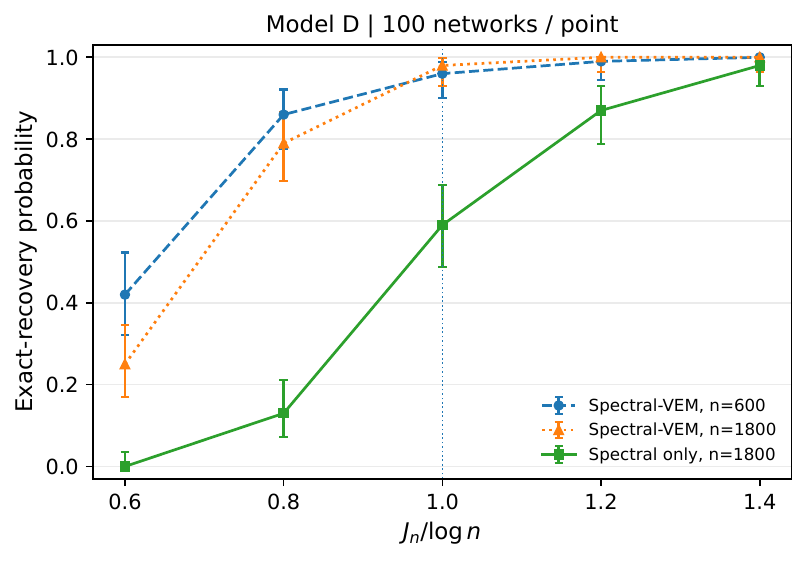}\hfill
\includegraphics[width=.485\textwidth]{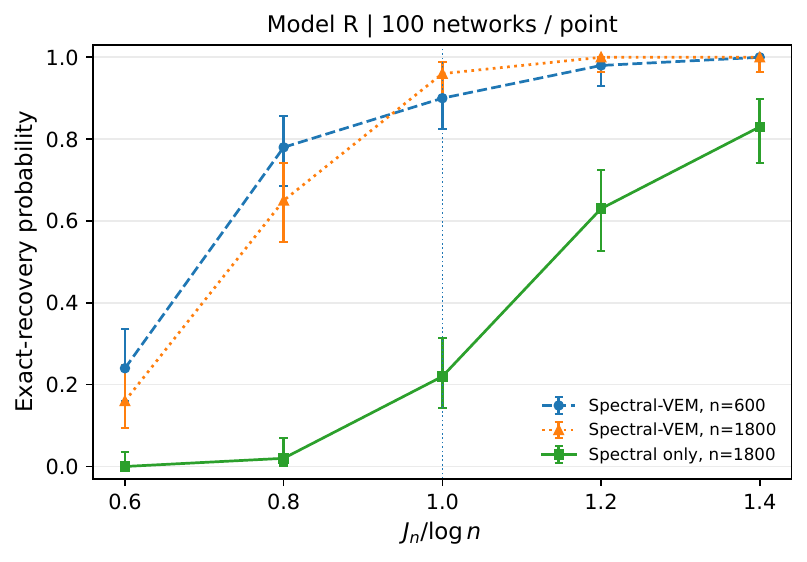}
\caption{Graph-level exact recovery near $J_n=\log n$. Each point is based on
100 independent graphs; intervals are two-sided exact 95\% binomial confidence
intervals. The dotted vertical line marks the first-order theoretical threshold
$J_n=\log n$. Additional counts are reported in the supplementary material.}
\label{fig:exact}
\end{figure}

Figure~\ref{fig:exact} reports exact recovery probabilities near the
first-order threshold. At $n=1800$ and $J_n=\log n$, SVEM recovers 97, 94,
98, and 96 of the 100 graphs exactly for A, H, D, and R, respectively. At
$J_n=1.2\log n$, the corresponding counts are 99, 100, 100, and 100.
Exact recovery becomes more frequent as $J_n/\log n$ increases, with the
transition occurring near the first-order reference $J_n=\log n$ in these
finite samples.

\subsection{Imbalance and initialization sensitivity}

\begin{figure}[tbp]\centering
\includegraphics[width=.485\textwidth]{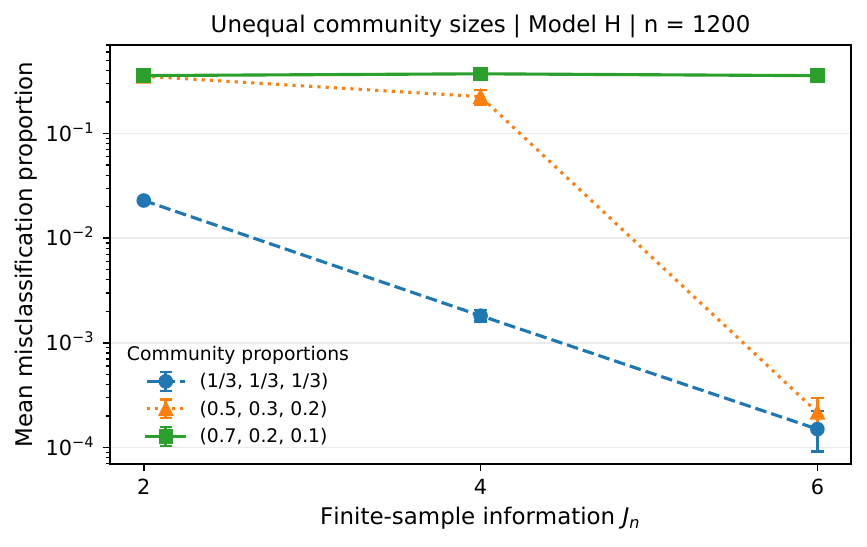}\hfill
\includegraphics[width=.485\textwidth]{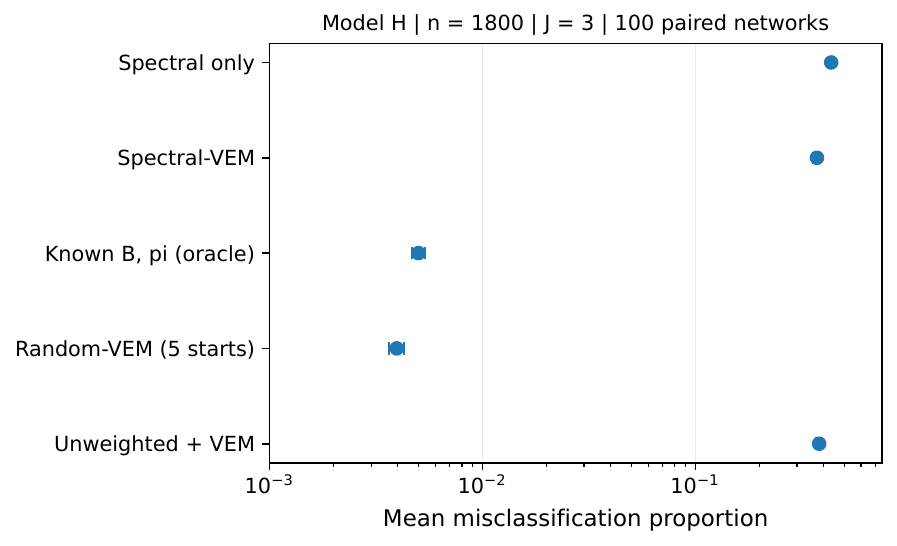}
\caption{Left: SVEM under three community-proportion vectors at $n=1200$.
Right: comparisons for H at $n=1800$ and $J_n=3$, using the same 100 graphs
for all procedures. Intervals are 95\% graph-level percentile bootstrap
intervals.}
\label{fig:imbalance}
\end{figure}

Matching $J_n$ does not remove finite sample sensitivity to community
imbalance. For H at $n=1200$ and $J_n=6$, the mean SVEM error is
$1.5\times10^{-4}$ for balanced proportions and $2.2\times10^{-4}$ for
$(0.5,0.3,0.2)$, but rises to $0.357$ for $(0.7,0.2,0.1)$. In the latter
setting, the spectral error is $0.433$, and none of the 100 runs satisfies
the stopping criterion. The left panel of
Figure~\ref{fig:imbalance} shows this sensitivity to community imbalance.

We next examine sensitivity to initialization and the spectral construction.
For H at $n=1800$ and $J_n=3$, the mean errors are $0.373$ for SVEM,
$0.0040$ for VEM selected from five random starts by the final working
objective, and $0.0050$ for the known-parameter reference. The large
difference indicates that, in this finite-sample setting, the prescribed
spectral initialization can place the iteration in an unfavorable region;
the random-start procedure also uses substantially more computation. The
known-parameter comparison shows that low error is attainable when $B$ and
$\pi$ are supplied, although it does not by itself identify the source of
the poor SVEM performance. For D at the same information level, using the
eigenvalues of largest magnitude gives mean error $0.0045$, compared with
$0.0327$ when the largest algebraic eigenvalues are used. For the rank-two
model R, weighted and unweighted spectral starts followed by VEM give the
same mean error in this setting. The right panel of Figure~\ref{fig:imbalance} 
gives the corresponding initialization comparisons for H.

To separate the effect of spectral initialization from the behavior of the
VEM update, we next vary the initial error directly while keeping the graph
and the update rule fixed.

\begin{figure}[tbp]\centering
\includegraphics[width=.485\textwidth]{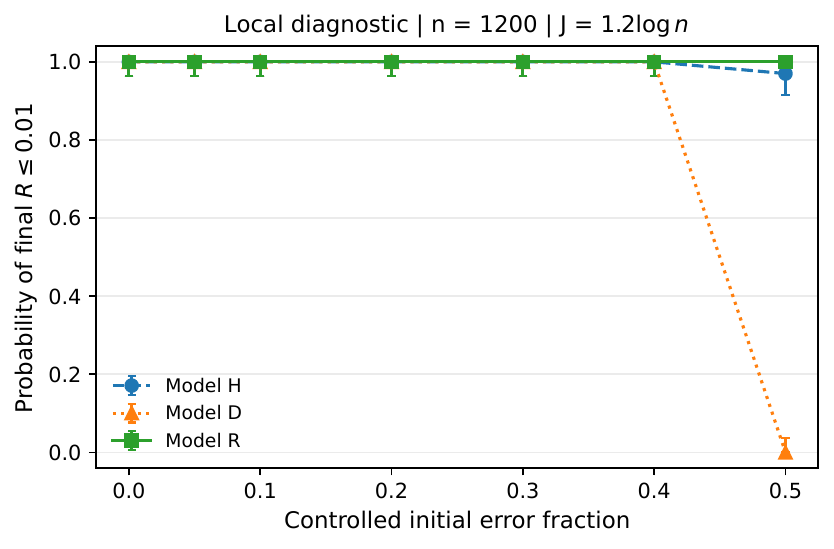}\hfill
\includegraphics[width=.485\textwidth]{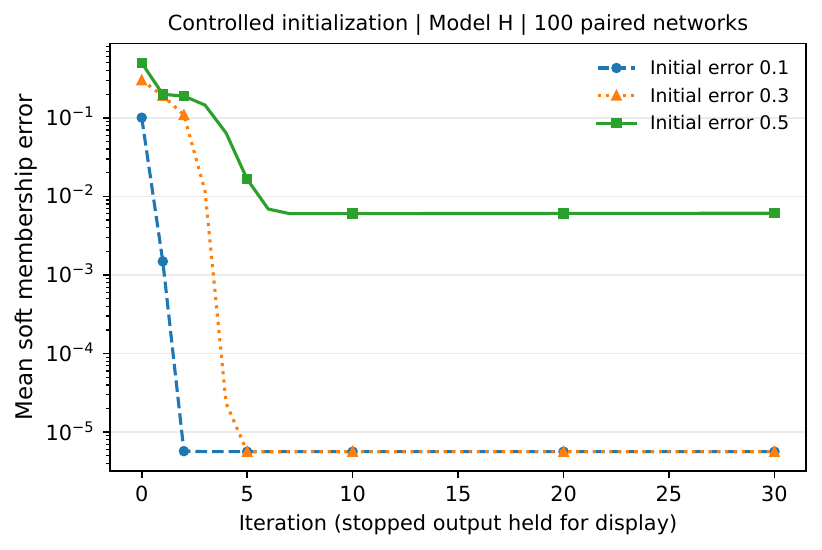}
\caption{Controlled initialization at $n=1200$ and $J_n=1.2\log n$, using
100 independent graphs per model. Left: probability of final $R\le0.01$
with exact binomial confidence intervals. Right: mean soft membership error
across iterations for H.}
\label{fig:basin}
\end{figure}

At $J_n=1.2\log(1200)$, H recovers all 100 graphs exactly for each tested
$\alpha\le0.4$, while D and R each recover 97 of 100. At $\alpha=0.5$, D
fails on all 100 graphs, whereas H and R each recover 97 of 100. The left
panel of Figure~\ref{fig:basin} reports the probability of final
$R\le0.01$ across the controlled initial-error levels. The right panel shows
the mean soft membership error across iterations for H. The controlled
initializations are constructed from the true labels solely for this
sensitivity analysis.

\subsection{Stopping behavior and computation}

Among the 8,500 primary SVEM runs, 1,527 reach the iteration limit without
satisfying the stopping criterion. None of the 6,973 runs that satisfy the
criterion changes its hard labels during the additional 20 updates, whereas
767 of the runs that reach the iteration limit do. Thus the tolerance based
stopping rule gives stable outputs when it is satisfied, while runs that reach
the iteration limit may still change with further iteration.

We also compare the practical stopping rule with the theoretical iteration
budget in Algorithm~\ref{alg:vem}. For the 1,600 graphs in the $n=600$ risk
settings, 1,167 runs satisfy the practical stopping criterion; none of their
misclassification rates changes when the iteration is continued to the
prescribed budget. Among the remaining 433 runs, the misclassification rate
decreases in 154, increases in 49, and is unchanged in 230. Detailed results
are given in Section~K of the supplementary material.

Computation times are reported in Table~\ref{tab:timing}. The median
time per VEM update increases from approximately $0.17$ milliseconds (ms) at $n=600$ to
$1.30$ ms at $n=4800$ as the network size and number of edges increase.
No vertex was removed by the trimming step in these simulations, so its
finite-sample effect is not examined here. Additional implementation and
reproducibility details are provided in Section~K of the supplementary material.

\begin{table}[tbp]
\centering\small\setlength{\tabcolsep}{5pt}
\begin{tabular}{rrrrrrr}\toprule
$n$ & Mean $|E|$ & Spectral (ms) & VEM (ms) & Total (ms)
& Per update (ms) & Iterations\\\midrule
600 & 15292 & 6.80 & 2.02 & 8.82 & 0.168 & 12 \\
1200 & 33982 & 9.08 & 3.83 & 12.91 & 0.319 & 12 \\
2400 & 74662 & 14.12 & 7.61 & 21.65 & 0.634 & 12 \\
4800 & 162664 & 24.44 & 15.61 & 40.07 & 1.301 & 12 \\
\bottomrule
\end{tabular}

\vspace{2mm}
\caption{Computation times for model H with
$\rho_n=8\log(n)/n$, based on 100 graphs at each network size.
Reported times and iteration counts are medians. Total time includes spectral
initialization and VEM but excludes the additional
20 updates used to assess the stopping criterion.}
\label{tab:timing}
\end{table}

\section{Discussion}\label{sec:discussion}

The main conclusion is that the prescribed spectrally initialized batch VEM
procedure can attain the oracle Chernoff exponent while estimating the full
block matrix and community proportions from the same graph. These guarantees
concern the estimator returned after the prescribed iteration budget and do
not assume global maximization of the variational objective. Nor do they
require monotone ascent of the objective or convergence of the empirical map
to a unique fixed point. Relative to planted partition iteration theory,
the result allows a general connectivity matrix; relative to likelihood ratio
refinement based on graph splitting, it applies to the simultaneous soft
updates used by VEM. The main step in the analysis is a uniform one-step error
bound with a common random remainder whose expectation is controlled at the
Chernoff scale.

To assess optimality, we compare the risk of the prescribed procedure with the
minimax risk over local parameter spaces in which both the connectivity matrix
and the community counts are allowed to vary. On these spaces, the minimax
expected misclassification risk has the same first order logarithmic exponent
as the weighted Chernoff information, and the prescribed procedure attains this
exponent throughout the fixed profile sparse regime with diverging degree. 
The procedure also attains the known first order threshold for exact recovery.

The numerical results show that VEM can substantially reduce spectral error
when the spectral initialization is sufficiently accurate, including for
disassortative models and models with rank deficiency. At lower information
levels, however, the outcome can be sensitive to the initialization, the
spectral construction, and highly unequal community proportions. The exact
recovery experiments also show the expected increase in success probability
as the information level approaches and exceeds the first order threshold.

\begin{acks}[Acknowledgments]
We acknowledge the use of a large language model (LLM) for language editing
and improving the clarity of the manuscript. All AI-assisted revisions were
carefully reviewed for accuracy and originality. The authors take full
responsibility for the content of this work.

\vspace{2mm}
\noindent
The authors are listed in
alphabetical order. Jiangzhou Wang is the corresponding author.
\end{acks}

\begin{funding}
Wang's work was supported by the National Natural Science Foundation of China
(Grant Nos. 12671329 and 12301355); the Basic and Applied Basic Research
Foundation of Guangdong Province (Grant No. 2024A1515011510); the China
Postdoctoral Science Foundation (Grant Nos. 2022T150293 and 2021M701588); and
the Shenzhen Science and Technology Program (Grant No. QNXMB20250701094201002).
\end{funding}

\begin{supplement}
\stitle{Supplement to ``Optimal Community Recovery by Spectrally Initialized
Variational EM in General Stochastic Block Models''}
\sdescription{Complete proofs of the theoretical results, additional numerical
results, and implementation details.}
\end{supplement}


\bibliographystyle{spr-ims-nameyear}

\bibliography{references}
\end{document}